\documentclass[12pt]{amsart}
\usepackage{amsmath,amssymb,amsbsy,amsfonts,amsthm,latexsym,amsopn,amstext,
                                               amsxtra,euscript,amscd,color,bm}
                                                                                           
\usepackage[colorlinks,linkcolor=blue,anchorcolor=blue,citecolor=blue,backref=page]{hyperref}

\hypersetup{breaklinks=true}

\begin{document}

\def\ov#1{{\overline{#1}}} 
\def\un#1{{\underline{#1}}}
\def\wh#1{{\widehat{#1}}}
\def\wt#1{{\widetilde{#1}}}

\newcommand{\Ch}{{\operatorname{Ch}}}
\newcommand{\Elim}{{\operatorname{Elim}}}
\newcommand{\proj}{{\operatorname{proj}}}
\newcommand{\h}{{\operatorname{h}}}

\newcommand{\hh}{\mathrm{h}}
\newcommand{\aff}{\mathrm{aff}}
\newcommand{\Spec}{{\operatorname{Spec}}}
\newcommand{\Res}{{\operatorname{Res}}}
\newcommand{\Orb}{{\operatorname{Orb}}}

\renewcommand*{\backref}[1]{}
\renewcommand*{\backrefalt}[4]{%
    \ifcase #1 (Not cited.)%
    \or        (p.\,#2)%
    \else      (pp.\,#2)%
    \fi}
\def\lc{{\mathrm{lc}}}
\newcommand{\hcan}{{\operatorname{\wh h}}}

\newcommand{\hooklongrightarrow}{\lhook\joinrel\longrightarrow}

\newcommand{\bfa}{{\boldsymbol{a}}}
\newcommand{\bfb}{{\boldsymbol{b}}}
\newcommand{\bfc}{{\boldsymbol{c}}}
\newcommand{\bfd}{{\boldsymbol{d}}}
\newcommand{\bff}{{\boldsymbol{f}}}
\newcommand{\bfg}{{\boldsymbol{g}}}
\newcommand{\bfell}{{\boldsymbol{\ell}}}
\newcommand{\bfp}{{\boldsymbol{p}}}
\newcommand{\bfq}{{\boldsymbol{q}}}
\newcommand{\bfs}{{\boldsymbol{s}}}
\newcommand{\bft}{{\boldsymbol{t}}}
\newcommand{\bfu}{{\boldsymbol{u}}}
\newcommand{\bfv}{{\boldsymbol{v}}}
\newcommand{\bfw}{{\boldsymbol{w}}}
\newcommand{\bfx}{{\boldsymbol{x}}}
\newcommand{\bfy}{{\boldsymbol{y}}}
\newcommand{\bfz}{{\boldsymbol{z}}}

\newcommand{\bfA}{{\boldsymbol{A}}}
\newcommand{\bfF}{{\boldsymbol{F}}}
\newcommand{\bfG}{{\boldsymbol{G}}}
\newcommand{\bfR}{{\boldsymbol{R}}}
\newcommand{\bfQ}{{\boldsymbol{Q}}}
\newcommand{\bfT}{{\boldsymbol{T}}}
\newcommand{\bfU}{{\boldsymbol{U}}}
\newcommand{\bfX}{{\boldsymbol{X}}}
\newcommand{\bfY}{{\boldsymbol{Y}}}
\newcommand{\bfZ}{{\boldsymbol{Z}}}

\newcommand{\bfeta}{{\boldsymbol{\eta}}}
\newcommand{\bfxi}{{\boldsymbol{\xi}}}
\newcommand{\bfrho}{{\boldsymbol{\rho}}}

\newcommand{\PreP}{{\mathrm{PrePer}}}

\def\fM{{\mathfrak M}}
\def\fE{{\mathfrak E}}
\def\fF{{\mathfrak F}}



\newfont{\teneufm}{eufm10}
\newfont{\seveneufm}{eufm7}
\newfont{\fiveeufm}{eufm5}
%
%
\newfam\eufmfam
                \textfont\eufmfam=\teneufm \scriptfont\eufmfam=\seveneufm
                \scriptscriptfont\eufmfam=\fiveeufm
%
%
\def\frak#1{{\fam\eufmfam\relax#1}}
%

\def\ts{\thinspace}

\newtheorem{theorem}{Theorem}
\newtheorem{lemma}[theorem]{Lemma}
\newtheorem{claim}[theorem]{Claim}
\newtheorem{cor}[theorem]{Corollary}
\newtheorem{prop}[theorem]{Proposition}
\newtheorem{question}[theorem]{Open Question}

\newtheorem{rem}[theorem]{Remark}
\newtheorem{definition}[theorem]{Definition}

\newenvironment{dedication}
  {
   \thispagestyle{empty}
   \vspace*{\stretch{1}}
   \itshape             
   \raggedleft          
  }
  {\par 
   \vspace{\stretch{3}} 
  }

\numberwithin{table}{section}
\numberwithin{equation}{section}
\numberwithin{figure}{section}
\numberwithin{theorem}{section}


\def\squareforqed{\hbox{\rlap{$\sqcap$}$\sqcup$}}
\def\qed{\ifmmode\squareforqed\else{\unskip\nobreak\hfil
\penalty50\hskip1em\null\nobreak\hfil\squareforqed
\parfillskip=0pt\finalhyphendemerits=0\endgraf}\fi}

\def\fA{{\mathfrak A}}
\def\fB{{\mathfrak B}}

\def\cA{{\mathcal A}}
\def\cB{{\mathcal B}}
\def\cC{{\mathcal C}}
\def\cD{{\mathcal D}}
\def\cE{{\mathcal E}}
\def\cF{{\mathcal F}}
\def\cG{{\mathcal G}}
\def\cH{{\mathcal H}}
\def\cI{{\mathcal I}}
\def\cJ{{\mathcal J}}
\def\cK{{\mathcal K}}
\def\cL{{\mathcal L}}
\def\cM{{\mathcal M}}
\def\cN{{\mathcal N}}
\def\cO{{\mathcal O}}
\def\cP{{\mathcal P}}
\def\cQ{{\mathcal Q}}
\def\cR{{\mathcal R}}
\def\cS{{\mathcal S}}
\def\cT{{\mathcal T}}
\def\cU{{\mathcal U}}
\def\cV{{\mathcal V}}
\def\cW{{\mathcal W}}
\def\cX{{\mathcal X}}
\def\cY{{\mathcal Y}}
\def\cZ{{\mathcal Z}}

\def\nrp#1{\left\|#1\right\|_p}
\def\nrq#1{\left\|#1\right\|_m}
\def\nrqk#1{\left\|#1\right\|_{m_k}}
\def\Ln#1{\mbox{\rm {Ln}}\,#1}
\def\nd{\hspace{-1.2mm}}
\def\ord{{\mathrm{ord}}}
\def\Cc{{\mathrm C}}
\def\Pb{\,{\mathbf P}}

\def\va{{\mathbf{a}}}

\newcommand{\commF}[1]{\marginpar{%
\begin{color}{red}
\vskip-\baselineskip 
\raggedright\footnotesize
\itshape\hrule \smallskip F: #1\par\smallskip\hrule\end{color}}}

\newcommand{\commI}[1]{\marginpar{%
\begin{color}{blue}
\vskip-\baselineskip 
\raggedright\footnotesize
\itshape\hrule \smallskip I: #1\par\smallskip\hrule\end{color}}}




\newcommand{\ignore}[1]{}

\def\vec#1{\boldsymbol{#1}}

\def\e{\mathbf{e}}



\def\GL{\mathrm{GL}}

\hyphenation{re-pub-lished}

\def\rank{{\mathrm{rk}\,}}
\def\dd{{\mathrm{dyndeg}\,}}
\def\lcm{{\mathrm{lcm}\,}}

\def\A{\mathbb{A}}
\def\B{\mathbf{B}}
\def \C{\mathbb{C}}
\def \F{\mathbb{F}}
\def \K{\mathbb{K}}
\def \L{\mathbb{L}}
\def \Z{\mathbb{Z}}
\def \P{\mathbb{P}}
\def \R{\mathbb{R}}
\def \Q{\mathbb{Q}}
\def \N{\mathbb{N}}
\def \U{\mathbb{U}}
\def \Z{\mathbb{Z}}

\def \nd{{\, | \hspace{-1.5 mm}/\,}}

\def\mand{\qquad\mbox{and}\qquad}

\def\Zn{\Z_n}

\def\Fp{\F_p}
\def\Fq{\F_q}
\def \fp{\Fp^*}
\def\\{\cr}
\def\({\left(}
\def\){\right)}
\def\fl#1{\left\lfloor#1\right\rfloor}
\def\rf#1{\left\lceil#1\right\rceil}
\def\vh{\mathbf{h}}
\def\ov#1{{\overline{#1}}}
\def\un#1{{\underline{#1}}}
\def\wh#1{{\widehat{#1}}}
\def\wt#1{{\widetilde{#1}}}
\newcommand{\abs}[1]{\left| #1 \right|}

\def\ZK{\Z_\K}
\def\LH{\cL_H}

\def \fI{\mathfrak{I}}
\def \fJ{\mathfrak{J}}
\def \fV{\mathfrak{V}}

\title[Counting Lattices]
{On the error term of a lattice counting problem}

\author{Florian~Luca}

\address{School of Mathematics, University of the Witwatersrand, Private Bag X3, Wits 2050, South Africa; Max Planck Institute for Mathematics, Vivatsgasse 7, 53111 Bonn, Germany; 
Department of Mathematics, Faculty of Sciences, University of Ostrava, 30 Dubna 22, 701 03
Ostrava 1, Czech Republic}
\email{florian.luca@wits.ac.za}

\author{Igor E. Shparlinski}

\address{Department of Pure Mathematics, University of New South Wales\\
2052 NSW, Australia.}
\email{igor.shparlinski@unsw.edu.au}



\date{}

\begin{abstract} We improve the error terms of some estimates related to counting lattices from  recent work of L.~Fukshansky,  P.~Guerzhoy  and F.~Luca   (2017).
This improvement is based on some analytic techniques, in particular on bounds
of exponential sums coupled with the use of Vaaler polynomials. 
\end{abstract}

\maketitle

\section{Introduction}

\subsection{Background}

For integer $T \ge 1$, we let
$$
\cF(T) = \{a/b~:~(a,b)\in \Z^2, \ 0 \le a< b\le T, \ \gcd(a,b)=1\}
$$ 
be the set of Farey fractions. 

Now, following~\cite{FGL}, we consider the quantity
$$
C(T) = \sum_{a/b \in \cF(T) \cap[0,1/2]} \#  \cC_{a,b}(T),
$$
where 
$$
\cC_{a,b}(T) =  \cF(T) \cap [1-a^2/b^2,1].
 $$

The quantity  $C(T)$ appears naturally in some counting 
 problems for two-dimensional lattices. More precisely, every similarity class of {\it planar lattices} 
 can be parametrised by a point $\tau=x_0+iy_0$ in 
 $$
\cR=\{\tau=x_0+i y_0 ~:~\quad  0\le x_0\le 1/2,\quad y_0\ge 0,\quad |\tau|\ge 1\} \subseteq \C, 
 $$
 where one identifies $\tau\in {\mathcal R}$ with the lattice
 $$
 \Lambda_{\tau}=\left(\begin{matrix} 1 & x_0\\ 0 & y_0\end{matrix}\right) {\mathbb Z}^2.
 $$
 Further, similarity classes of {\it arithmetic} planar lattices correspond to $\Lambda_{\tau}$, where 
 $$
 \tau 
 =a/b+i \sqrt{c/d}
 $$
 for integers $a,b,c,d$ such that
 $$
 \gcd(a,b)=\gcd(c,d)=1, \quad 0\le a\le b/2,\quad d>0, \quad c/d\ge 1-a^2/b^2.
 $$
The class is {\it semistable} if furthermore $c\le d$. 
 With these conventions, the quantity $C(T)$ counts the number of similarity classes of semi-stable arithmetic planar lattices of height at most $T$, 
 that if for which $\max\{a,b,c,d\}\le T$.  
 
The following result appears as~\cite[Lemma~3.2]{FGL}: 
\begin{equation}
\label{eq:FGL}
C(T)  =  \frac{3}{8\pi^4} T^4 + O(T^3 \log T).
\end{equation} 

Our goal here is to sharpen the error term in the asymptotic formula~\eqref{eq:FGL} 
and in particular we show that that error term can be taken to be 
$O\(T^3(\log T)^{2/3} (\log\log T)^{1+o(1)}\)$ (see Corollary~\ref{cor:CT expl} below). 
However, it seems to be more natural 
to express the main term via some general quantities related to Farey fractions 
and then try to minimize the error term.  In particular, we outline some results
on counting Farey fractions in Section~\ref{sec:Farey}. 

Here, we accept this point of view and 
thus express the main term of the asymptotic formula for $\#C(T)$ 
via the cardinality 
$$
F(T) = \# \cF(T)
$$
of the set of of Farey fractions and also second moment 
of the Farey fractions in $[0,1/2]$:
$$
G(T) = \sum_{\substack{\xi \in \cF(T)\\\xi\le1/2}} \xi^2, \qquad \nu = 0, 1, \ldots.
$$

It is also convenient to define 
\begin{equation}
\label{eq:Mt}
 M(t) = \sum_{1\le k \le t} \mu(k).
 \end{equation}

As usual $A= O(B)$,  $A \ll B$, $B \gg A$ are equivalent to $|A| \le c |B|$ for some 
 {\it absolute\/}  constant $c> 0$, whereas $A=o(B)$ means that $A/B\to 0$.

\begin{theorem}
\label{thm:CT FTGT}
We have 
$$
C(T)  =  F(T)G(T) + O\left(T^{11/4+o(1)}+T^3\delta(T^{1/2})\log T\right), 
 $$
where $\delta(t)$ is any  decreasing function such that 
$$
M(t) \le t \delta(t)
$$
holds. 
\end{theorem}

By the classical bound of Walfisz~\cite[Chapter~V, Section~5, Equation~(12)]{Wa} one can take
\begin{equation}
\label{eq:MobSum Uncond}
 \delta(t) = \exp(-c(\log T)^{3/5} (\log\log T)^{-1/5})
 \end{equation} 
 for absolute constant $c>0$, 
hence immediately producing the bound $O(T^3 \exp(- c_0(\log T)^{3/5} (\log\log T)^{-1/5}))$
for some constant $c_0>0$
on  the error term in Theorem~\ref{thm:CT FTGT}. 
Under the Riemann Hypothesis, we can take 
\begin{equation}
\label{eq:MobSum RH} 
 \delta(t) = t^{-1/2+\varepsilon}
 \end{equation} 
for any $\varepsilon>0$ (see~\cite{Titch}). Without $\varepsilon$, the inequality \eqref{eq:MobSum RH}  is known as a conjecture of Mertens which has been refuted 
by Odlyzko and te Riele~\cite{Odtr}. Hence, under the Riemann Hypothesis we 
obtain an error  
$O(T^{11/4+o(1)})$ as $T\to\infty$.   In~\eqref{eq:GTFT} below we obtain an approximation to $G(T)$
via $F(T)$
which implies the following result.

\begin{cor}
\label{cor:CT FT}
We have 
$$
C(T)  = \frac{1}{24} F(T)^2 + O(T^3).
 $$
\end{cor}

Finally, using the asymptotic formula for $F(T)$ with the error term 
given by~\eqref{eq:Walf}, we obtain the following direct improvement of~\eqref{eq:FGL}: 

\begin{cor}
\label{cor:CT expl}
We have 
$$
C(T)  = \frac{3}{8\pi^4} T^4 + O(T^3(\log T)^{2/3} (\log\log T)^{1+o(1)})
 $$
 as $T\to\infty$.
\end{cor}

We remark that improving the error term in  Corollary~\ref{cor:CT expl}
is probably impossible until the bound~\eqref{eq:Walf} is improved. 
However, it is plausible that one can improve~\eqref{eq:GTFT} and 
thus obtain a stronger version of Corollary~\ref{cor:CT FT}, which we pose
as an open question. 
 
\section{Main Term}

\subsection{Initial transformations}

By a result of Niederreiter~\cite{Nied1}, for any integers $0 \le a <b$ the following formula holds
\begin{equation}
\label{eq:Nied1}
 \# \cC_{a,b}(T)  - \frac{a^2}{b^2}  F(T)  =  - \sum_{n=1}^T \sum_{d\mid n} \mu(n/d) \{da^2/b^2\},
\end{equation} 
where $\mu(k)$ is the M{\"o}bius function   (see~\cite[Equation~(1.16]{IwKow})
and $\{\alpha\}$ is the fractional 
part of a real $\alpha$. 

We rewrite~\eqref{eq:Nied1} as
$$
 \# \cC_{a,b}(T)  - \frac{a^2}{b^2}  F(T)  =  - \sum_{d=1}^T \{da^2/b^2\} M(T/d) ,
$$
where $M(t)$ is given by~\eqref{eq:Mt}.
We now write
 \begin{equation}
\label{eq:Nied3}
C(T)  - \fM(T)  = \fE(T),
\end{equation} 
where
\begin{align*}
 \fM(T) & =     F(T)  \sum_{a/b \in \cF(T) \cap[0,1/2]} \frac{a^2}{b^2} = F(T) G(T),  \\
\fE(T) & =  -  \sum_{a/b \in \cF(T) \cap[0,1/2]}   \sum_{d=1}^T \{db^2/b^2\} M(T/d) .
\end{align*} 
Using either of the bounds~\eqref{eq:MobSum Uncond} and~\eqref{eq:MobSum RH}
gives the bound $O(T)$ for each inner sum in the definition of the error term $\fE(T)$
(see, for example, the proof of~\cite[Lemma~2]{Nied1}), 
and thus 
yields the conclusion of Theorem~\ref{thm:CT FTGT} 
with an error term  $O(T^3)$. Thus,  to do better, we need to 
investigate the {\it cancellations\/} between these sums.

\subsection{Counting Farey fractions}
\label{sec:Farey}

Here, we collect some known facts about  Farey fractions.

The set $\cF(T)$ has been the subject of a lot of research. Writing $\varphi(n)$ for the Euler function of the positive integer $n$, we have
$$
F(T) =\sum_{b\le T} \varphi(b)=\frac{3}{\pi^2} T^2+R(T).
$$
The error term $R(T)$ above has also been the subject of a lot of research. For example, 
by the classical result of Mertens~\cite{Me}  (that dates back to~1874), we have
$$
R(T)=O(T\log T). 
$$ 
This has been improved by Walfisz~\cite[Chapter~V, Section~5, Equation~(35)]{Wa} and then finally by Saltykov~\cite{Sal}  to 
 \begin{equation}
\label{eq:Walf}
R(T)=O\(T(\log T)^{2/3} (\log\log T)^{1+o(1)}\)
\end{equation} 
as $T\to \infty$.  

Erd\H os and Shapiro~\cite{ES} have shown that 
 $$
 R(T)=\Omega_{\pm }(T\log\log\log\log T),
 $$ 
 which means that for some positive constant $c$, each of the inequalities
$$
R(T)>c T\log\log\log\log T\qquad {\text{\rm and}}\qquad R(T)<-cT\log\log\log\log T
$$
holds infinitely often, while Montgomery~\cite{Mo} has sharpened this to 
$$
R(T)=\Omega_{\pm }(T (\log\log T)^{1/2}).
$$ 
Average values and moments of $R(T)$ have also been considered. For example,
\begin{equation}
\label{eq:1}
\sum_{m\le T} R(m)=\frac{3 T^2}{2\pi^2}+O(T^2\eta(T))
\end{equation}
(see~\cite{SuSi}), 
and 
\begin{equation}
\label{eq:2}
\sum_{m\le T} R(m)^2=\left(\frac{1}{6\pi^2}+\frac{2}{\pi^4}\right)T^4+O(T^3\eta(T)),
\end{equation}
(see~\cite{BRL}), where in both~\eqref{eq:1} and~\eqref{eq:2} 
$$
\eta(T)=\exp(-A(\log T)^{3/5}(\log\log T)^{-1/5})
$$
for some constant $A>0$ (not necessarily the same one in both~\eqref{eq:1} and~\eqref{eq:2}). 

We remark that for the second (and 
other) moments of Farey fractions one can obtain asymptotic formulas via the general bounds on the difference between 
sums of continuous functions on Farey fractions and the corresponding integrals (see~\cite{BKY,Yos}). 

Unfortunately, these results do not seem to apply to the sum 
$G(T)$. On the other hand,  one can, via elementary but rather tedious arguments, 
relate $G(T)$ to $F(T)$ and then show 
that 
\begin{equation}
\label{eq:GT asymp}
G(T) =
 \frac{1}{8\pi^2}  T^2 + O\(T(\log T)^{2/3} (\log\log T)^{1+o(1)}\)
\end{equation}
as $T\to\infty$. However, here we use some general results to derive~\eqref{eq:GT asymp}.
We start with recalling  the bound 
$$
\Delta(T) = O(T^{-1})
$$ 
of Niederreiter~\cite{Nied1} on the discrepancy
$$
\Delta(T)  = \sup_{0 \le \alpha \le 1}
\left| \#\(\cF(T) \cap [0,\alpha]\) - \alpha F(T)\right|
$$
of the Farey fractions. 

Since the function 
$$
f(z) =   \left\{ \begin{array}{ll} z^2&   \text{ if}\ z\in [0,1/2], \\ 
0&   \text{ if}\ z\in (1/2,1],\end{array} \right.
$$
is of bounded variation, by the classical {\it Koksma inequality\/} 
(see, for example,~\cite[Theorem~2.9]{Nied2}), we have 
 \begin{equation}
 \begin{split}
\label{eq:GTFT}
G(T) & = \sum_{\xi \in \cF(T)} f(\xi) \\
& = F(T) \int_0^1 f(z) dz + O\(F(T)\Delta(T)\) =\frac{1}{24} F(T) + O(T),
\end{split}
\end{equation}
which together with~\eqref{eq:Walf} implies the estimate~\eqref{eq:GT asymp}. 

Finally, the asymptotic formulas~\eqref{eq:Walf} and~\eqref{eq:GT asymp}
imply Corollary~\ref{cor:CT expl}.

\section{Error Term}

\subsection{Some sums with  the M{\"o}bius function}

In handling the sums $\fM(T)$ and $\fE(T)$
we often appeal to a result of Gupta~\cite{Gupta}:

\begin{lemma}
\label{lem:Gupta}
For any integer $m\ge 1$, we have 
$$
\sum_{\substack{d=1\\\gcd(d,m)=1}}^T   \mu(d)\fl{T/d} =
\sum_{\substack{d \mid m^\ell\\ d \le T}} 1,
$$
where 
$$
\ell = \fl{\frac{\log T}{\log 2}}.
$$
\end{lemma}

Note that after changing the order summations,  Lemma~\ref{lem:Gupta} yields 
$$
\sum_{b \le T}
\sum_{\substack{d\mid b\\ \gcd(d,m)=1}} \mu(d) d
= 
\sum_{\substack{d=1\\\gcd(d,m)=1}}^T   \mu(d)\fl{T/d} =
\sum_{\substack{d \mid m^\ell\\ d \le T}} 1. 
$$
Thus, using it for $m=1$, we obtain:

\begin{cor}
\label{cor:sum div mu}
For the following sums we have 
$$
\sum_{b \le T} \sum_{d\mid b} \mu(d) d 
=1.
$$
\end{cor}

We remark, that somewhat related sums have also appeared in the work of Kunik~\cite{Kun1,Kun2}. 
However, these sums are independent  and thus our approach is different and in 
particular allows for a power saving, while the sums in~\cite{Kun2}  
are estimated with a much weaker saving.

\subsection{Vaaler polynomials}

We define the functions
$$\psi(u) = \{u\}-1/2 \mand \e(u) = \exp(2\pi i u). 
$$
By a result of Vaaler~\cite{Vaal} (see also~\cite[Theorem~A.6]{GrKol}), we have: 

\begin{lemma}
\label{lem:Vaal Approx}
For any integer $H\ge 1$ there is a trigonometric 
polynomial
$$
\psi_H(u)  = \sum_{1\leq \abs{h} \leq H} \frac{a_h}{-2i\pi h} \e(hu)
$$
with coefficients $a_h\in[0,1]$ and such 
that 
$$
\abs{\psi(u)-\psi_H(u)}
\le \frac{1}{2H+2} \sum_{\abs{h}\leq H} \(1-\frac{\abs{h}}{H+1}\) \e(hu).
$$
\end{lemma}

We now note that, by Lemma~\ref{lem:Vaal Approx}, we have
\begin{equation}
\label{eq:NiedVaal1}
 \sum_{d=1}^T \{da^2/b^2\} M(T/d) = \fE_0 + O(\fE_1 + \fE_2 + T^3/H), 
 \end{equation} 
where 
\begin{align*}
\fE_0 & =  \frac{1}{2}\sum_{a/b \in \cF(T) \cap[0,1/2]}  1 \times    \sum_{d=1}^T M(T/d),\\
 \fE_1 & =  \sum_{1\leq \abs{h} \leq H}  \abs{a_h}  
  \abs{  \sum_{a/b \in \cF(T) \cap[0,1/2]}   \sum_{d=1}^T  M(T/d)
 \e(a^2dh/b^2)},\\
 \fE_2 & = H^{-1} \sum_{1\leq \abs{h} \leq H}   
  \abs{  \sum_{a/b \in \cF(T) \cap[0,1/2]} \sum_{d=1}^T  M(T/d)
 \e(a^2dh/b^2)}\\
\end{align*}
(note that $ T^3/H$ comes from the contribution of the term with $h=0$ on
the right hand side of the inequality of Lemma~\ref{lem:Vaal Approx}). 

Clearly,
 $$
\fE_0 =  - \(\frac{1}{4} \cF(T) + O(1)\)   \sum_{d=1}^T M(T/d).
 $$
 Rearranging, for every  integer $T \ge 1$, we obtain 
 $$
 \sum_{d=1}^T M(T/d) =  \sum_{k=1}^T   \mu(k)\fl{T/k} = 1,
$$
by Corollary~\ref{cor:sum div mu}. 
Hence, 
  \begin{equation}
\label{eq:BoundM}
\fE_0 \ll T^2. 
\end{equation}

Substituting~\eqref{eq:BoundM}  in~\eqref{eq:NiedVaal1}
and combining this with~\eqref{eq:Nied3} we obtain 
\begin{equation}
\label{eq:NiedVaal2}
\fE(T)
  \ll \fE_1 + \fE_2 + T^3/H + T^2, 
 \end{equation}

\subsection{Bounds of exponential sums}

 Let 
 $$
 J = \fl{\frac{\log T}{\log 2}}.
 $$
We also fix two more positive integer  parameters $H\le T$ and $I \le J$, 
to be determined later.

We  fix some  parameters
 Define 
 $$
 \cD_i =   \Z \cap\left[2^i, \max\{T,2^{i+1}\}\right], \qquad i = I, \ldots, J.
 $$
 Using the definition of $\delta(t)$, we have
 \begin{equation}
 \begin{split}
\label{eq:E2F}
 & \fE_1 \ll    \sum_{1\leq \abs{h} \leq H} \frac{1}{h}
 \sum_{i=I}^J  |W_{h,i} |  + T^3  \delta(T/2^I)  \log H,\\
 & \fE_2\ll \frac{1}{H} \sum_{1\leq \abs{h} \leq H}  
 \sum_{i=I}^J  |W_{h,i}|   + T^3  \delta(T/2^I) , 
\end{split}
\end{equation} 
where
$$
W_{h,i} =   \sum_{a/b \in \cF(T) \cap[0,1/2]}   \sum_{d \in \cD_i}  M(T/d)
 \e(a^2dh/b^2), \qquad i =I, \ldots, J.
$$
We  fix $i \in [I,J]$ and write 
$$
W_{h,i} =  \sum_{d \in \cD_i}  M(T/d)  \sum_{b=1}^T
\sum_{\substack{1 \le a \le b/2\\ \gcd(a,b)=1} }    \e(a^2dh/b^2).
$$
We estimate $M(T/d)$ trivially as
$$
|M(T/d)| \le T/d \ll T2^{-i},
$$
and obtain 
$$
W_{h,i} =  T2^{-i} \sum_{d \in \cD_i} \sum_{b=1}^T \left| 
\sum_{\substack{1 \le a \le b/2\\ \gcd(a,b)=1} }    \e(a^2dh/b^2)\right|.
$$
Using that $ \#\cD_i \ll 2^i$, by the Cauchy inequality,  we obtain
$$
|W_{h,i}|^2  \ll   T^3 2^{-i}  \sum_{d \in \cD_i}^T 
\sum_{b=1}^T \left|
\sum_{\substack{1 \le a \le b/2\\ \gcd(a,b)=1} }    \e(a^2dh/b^2)\right|^2.
$$
Squaring out and changing the order of summations yields 
$$
|W_{h,i}|^2  \ll   T^3 2^{-i}  \sum_{b=1}^T
\sum_{\substack{1 \le a,c \le b/2\\ \gcd(ac,b)=1} }   \sum_{d \in \cD_i}   \e((a^2-c^2)dh/b^2).
$$

For  integer $q$ and $u$ define
$$
\langle u\rangle_q =\| u- q\Z\| = \min_{k \in \Z} |u - kq|
$$
as the distance to the closest integer 
 which is a multiple of $q$.
Then 
$$
 \sum_{d \in \cD_i}   \e((a^2-c^2)dh/b^2) \ll  \min\left\{2^i,  \frac{b^2}{\langle (a^2-c^2)h\rangle_{b^2}}\right \} 
$$
(see~\cite[Bound~(8.6)]{IwKow}). Thus,
\begin{align*}
|W_{h,i}|^2 & \ll  T^3 2^{-i}    \sum_{b=1}^T
\sum_{1 \le a,c \le b}    \min\left\{2^i,  \frac{b^2}{\langle (a^2-c^2)h\rangle_{b^2}}\right \}\\
&  \ll T^3 2^{-i}    \sum_{b=1}^T
\sum_{1 \le a,c \le b}    \min\left\{2^i,  \frac{b^2}{\langle (a^2-c^2)h\rangle_{b^2}}\right \},
\end{align*} 
where we have dropped the coprimality condition and extended the summation up to $b$
(only for the sake typographical simplicity). 

It is convenient to estimate separately the  contribution from  the diagonal $ a = c$, which leads to 
  \begin{equation}
\label{eq:BoundW1}
|W_{h,i}|^2 \ll  T^3 2^{-i}   \sum_{b=1}^T  
\sum_{1 \le a < c \le b}    \min\left\{2^i,  \frac{b^2}{\langle (a^2-c^2)h\rangle_{b^2}}\right \}
+ T^5.
\end{equation}

Now for every integer  $b \in [1, T]$ we define the set
$$
\cZ_0(b)=\left\{z\in \Z~:~|z|\le  2^{-i} b^2 \right \}.
$$
Furthermore, for $j =0,  \ldots, J $, we define  the sets
$$
\cZ_j(b)= \left\{z\in \Z \cap [-b^2/2,b^2/2]~:~2^{j-i}b^2<|z|\le 2^{j-i+1} b^2\right \}.
$$

Next, we fix  some $h$  in the interval $1 \le h \le H$ and define the sets:
\begin{align*}
\cA_{j}(b)&=\{(a,c) \in \Z^2~:~1 \le a<c \le b, \\
 & \qquad \qquad (a^2-c^2)h \equiv z \pmod {b^2}\ \text{for some}\ z\in\cZ_j  \}.
\end{align*} 
In particular,
  \begin{equation}
\label{eq:Sum ac}
\sum_{1 \le a<c \le b}    \min\left\{2^i,  \frac{b^2}{\langle (a^2-c^2)h\rangle_{b^2}}\right \}
\ll \sum_{j=0}^J 2^{i-j} \# \cA_{j}(b).
\end{equation} 
To estimate  $\# \cA_{j}(b)$ we note that for each $z$ the congruence 
$$
(a^2-c^2)h \equiv z \pmod {b^2}
$$ 
puts $ a^2-c^2$ in $\gcd(h,b^2)$
arithmetic progressions modulo $b^2$. Since 
$0 < c^2-a^2 < b^2$, each of these progressions, leads to an equation $c^2 -a^2 = k$
with some  positive  integer $k\le b^2\le T^2$. Using the classical bound on the 
divisor function $\tau(m)$ of the integer $m$  (see~\cite[Equation~(1.81]{IwKow}), 
we obtain 
\begin{align*}
 \#\cA_{j}(b) & \le \gcd(h,b^2) \# \cZ_j(b) \max\{\tau(k): k\le T^2\}\\
 &\le \gcd(h,b^2) \# cZ_j(b) T^{o(1)} \\
 &  = \gcd(h,b^2) \(2^{j-i}b^{2} +1\) T^{o(1)}\\
 & \le  \gcd(h,b^2) 2^{j-i}T^{2 + o(1)},
\end{align*}
as $T\to\infty$.  Using this in~\eqref{eq:Sum ac}, we obtain
\begin{align*}
\sum_{1 \le a<c \le b}    \min & \left\{2^i,  \frac{b^2}{\langle (a^2-c^2)h\rangle_{b^2}}\right \}\\
& \ll  J \gcd(h,b^2) T^{2+o(1)}\ll \gcd(h,b^2) T^{2+o(1)},
\end{align*}
where we ignored the $J$ factor because of the presence of the factor $T^{o(1)}$.

 With this notation, we infer from~\eqref{eq:BoundW1} that 
  \begin{equation}
\label{eq:BoundW3}
|W_{h,i}|^2 \ll  T^{5+o(1)}  2^{-i}  \sum_{b=1}^T \gcd(h,b^2) + T^{5}. 
\end{equation} 
Since obviously 
$$  \frac{1}{H} \sum_{1\leq \abs{h} \leq H}  
 \sum_{i=I}^J  |W_{h,i}|  \le  \sum_{1\leq \abs{h} \leq H}  \frac{1}{h}
 \sum_{i=I}^J  |W_{h,i} |, 
$$
we derive from~\eqref{eq:NiedVaal2} and ~\eqref{eq:E2F} (and absorbing the term $T^2$ into $T^3/H$ 
as $H \le T$), that
\begin{equation}
\label{eq:ESigma}
\fE(T)
  \ll 2^{-I/2} T^{5/2+o(1)} \Sigma  + J T^{5/2} \log H + T^3  \delta(T/2^I)   \log H + T^3/H,
\end{equation} 
where 
$$
\Sigma   =   \sum_{1\leq \abs{h} \leq H} \frac{1}{h}   \(\sum_{b=1}^T \gcd(h,b^2) \)^{1/2}. 
$$
Writing $h^{-1} = h^{-1/2} h^{-1/2} $ and using the Cauchy inequality, we obtain 
$$
\Sigma^2   \ll   \log H\sum_{1\leq \abs{h} \leq H} \frac{1}{h}  \sum_{b=1}^T \gcd(h,b^2).
$$
Furthermore, changing the order of summation and collecting together, 
for each divisor $d \mid b^2$,  the values $h$ with $\gcd(h,b^2) = d$, we obtain
\begin{align*}
\sum_{1\leq \abs{h} \leq H} \frac{1}{h}  \sum_{b=1}^T \gcd(h,b^2)
&=  \sum_{b=1}^T \sum_{1\leq \abs{h} \leq H} \frac{1}{h}  \gcd(h,b^2)\\
& \le  \sum_{b=1}^T\sum_{d\mid b^2}  d \sum_{1\leq \abs{k} \leq H/d} \frac{1}{dk} =  \sum_{b=1}^T\sum_{d\mid b^2}   \sum_{1\leq \abs{k} \leq H/d} \frac{1}{k} \\
& \le  \log H \sum_{b=1}^T \tau(b^2)\ll T(\log H) (\log T)^2.
\end{align*}
For the last estimate above, apply the main result of \cite{LuTo} to the function  $f(n)=\tau(n^2)$ which satisfies the conditions of that theorem with $k=3$. 
Substituting this in~\eqref{eq:ESigma},  we obtain
$$
\fE(T)
  \ll 2^{-I/2} T^{3+o(1)}  +  T^3  \delta(T/2^I)  \log H  +  T^{5/2} ( \log H)^2  .
$$
Choosing now $H=T^{1/2}$ and defining $I$ by the inequalities 
$$
2^{I-1} < T^{1/2} \le 2^I,
$$
we get the conclusion of Theorem \ref{thm:CT FTGT}.

\section*{Acknowledgements} 

F. L. was  
supported in part by grant CPRR160325161141 and an A-rated scientist award 
both from the NRF of 
South Africa and by grant no. 17-02804S of the Czech Granting Agency.  
I. E. S. was  
supported in part by ARC (Australia)  Grant DP140100118. 

Both authors would like to thank
the Max Planck Institute for Mathematics, Bonn, for the generous 
support and hospitality.

\end{document}